\documentclass[onecolumn]{article}
\usepackage{amsfonts}
\usepackage{amssymb}
\usepackage{amsmath}
\usepackage[numbers]{natbib}

\newtheorem{theorem}{Theorem}

\newtheorem{corollary}{Corollary}

\newtheorem{definition}{Definition}
\newtheorem{example}{Example}

\newtheorem{lemma}{Lemma}

\newtheorem{proposition}{Proposition}

\input{tcilatex}
\begin{document}

\title{Lyapunov Exponents of Free Operators}
\author{Vladislav Kargin \thanks{%
Department of Mathematics, Stanford University; kargin@stanford.edu}}
\date{August 25, 2008}
\maketitle

\begin{center}
\textbf{Abstract}
\end{center}

\begin{quotation}
Lyapunov exponents of a dynamical system are a useful tool to gauge the
stability and complexity of the system. This paper offers a definition of
Lyapunov exponents for a sequence of free linear operators. The definition
is based on the concept of the extended Fuglede-Kadison determinant. We
establish the existence of Lyapunov exponents, derive formulas for their
calculation, and show that Lyapunov exponents of free variables are additive
with respect to operator product. We illustrate these results using an
example of free operators whose singular values are distributed by the
Marchenko-Pastur law, and relate this example to C. M. Newman's
\textquotedblleft triangle\textquotedblright\ law for the distribution of
Lyapunov exponents of large random matrices with independent Gaussian
entries. As an interesting by-product of our results, we derive a relation
between the extended Fuglede-Kadison determinant and Voiculescu's
S-transform.
\end{quotation}

\section{Introduction}

Suppose that at each moment of time, $t_{i},$ a system is described by a
state function $\varphi \left( t_{i}\right) $ and evolves according to the
law $\varphi \left( t_{i+1}\right) =X_{i}\varphi \left( t_{i}\right) ,$
where $X_{i}$ is a sequence of linear operators. One can ask how small
changes in the initial position of the system are reflected in its long-term
behavior. If operators $X_{i}$ do not depend on time, $X_{i}=X,$ then the
long-term behavior depends to a large extent on the spectrum of the operator 
$X.$ If operators $X_{i}$ do depend on time but can be modelled as a
stationary stochastic process, then the long-term behavior of the system
depends to a large extent on so-called \emph{Lyapunov exponents} of the
process $X_{i}$.

The largest Lyapunov exponent of a sequence of random matrices was
investigated in a pioneering paper \citep{furstenberg_kesten60} by
Furstenberg and Kesten. This study was followed in \citep{oseledec68} by
Oseledec, who researched other Lyapunov exponents and finer aspects of the
asymptotic behavior of matrix products. These investigations were greatly
expanded and clarified by many other researchers. In particular, in %
\citep{ruelle82} Ruelle developed a theory of Lyapunov exponents for random
compact linear operators acting on a Hilbert space. Lyapunov exponents for
random $N\times N$ matrices when $N\rightarrow \infty $ were studied in %
\citep{cohen_newman84}, \citep{newman86a}, \citep{newman86b}, and %
\citep{isopi_newman92}$.$

The goal of this paper is to investigate how the concept of Lyapunov
exponents can be extended to the case of free linear operators. It was noted
recently by Voiculescu (\citep{voiculescu91}) that the theory of free
operators can be a natural asymptotic approximation for the theory of large
random matrices. Moreover, it was noted that certain difficult calculations
from the theory of large random matrices become significantly simpler if
similar calculations are performed using free operators. For this reason it
is interesting to study whether the concept of Lyapunov exponents is
extendable to free operators, and what methods for calculation of Lyapunov
exponents are available in this setting.

Free operators are not random in the traditional sense so the usual
definition of Lyapunov exponents cannot be applied directly. Our definition
of Lyapunov exponents is based on the observation that in the case of random
matrices, the sum of logarithms of the $k$ largest Lyapunov exponents equals
the rate at which a random $k$-dimensional volume element grows
asymptotically when we consecutively apply operators $X_{i}$.

In the case of free operators we employ the same idea. However, in this case
we have to clarify how to measure the change in the "$t$-dimensional volume
element" after we apply operators $X_{i}.$ It turns out that we can measure
this change by a suitable extension of the Fuglede-Kadison determinant.
Given this extension, the definition proceeds as follows: Take a subspace of
the Hilbert space, such that the corresponding projection is free from all $%
X_{i}$ and has the dimension $t$ relative to the given trace. Next, act on
this subspace by the sequence of free operators $X_{i}$. Apply the
determinant to measure how the \textquotedblleft volume
element\textquotedblright\ in this subspace changes under these linear
transformations. Use the asymptotic growth in the determinant to define the
Lyapunov exponent corresponding to the dimension $t$.

It turns out that the growth of the $t$-dimensional volume element is
exponential with a rate which is a function of the dimension $t.$ We call
this rate the \emph{integrated Lyapunov exponent}. It is an analogue of the
sum of the $k$ largest Lyapunov exponents in the finite-dimensional case.
The derivative of this function is called the \emph{marginal Lyapunov
exponent}. Its value at a point $t$ is an analogue of the $k$-th largest
Lyapunov exponent.

Next, we relate the marginal Lyapunov exponent $f_{X}\left( t\right) $ to
the Voiculescu $S$-transform of the random variable $X_{i}^{\ast }X_{i}.$
The relationship is very simple: 
\begin{equation}
f_{X}\left( t\right) =-\left( 1/2\right) \log \left[ S_{X^{\ast }X}\left(
-t\right) \right] .  \label{formula_marginal_Lyapunov_and_S}
\end{equation}

Using this formula, we prove that the marginal Lyapunov exponent is
decreasing in $t,$ and derive an expression for the largest Lyapunov
exponent. Formula (\ref{formula_marginal_Lyapunov_and_S}) also allows us to
prove the additivity of the marginal Lyapunov exponent with respect to
operator product: If $X$ and $Y$ are free, then $f_{XY}\left( t\right)
=f_{X}\left( t\right) +f_{Y}\left( t\right) .$

As an example of application of formula (\ref%
{formula_marginal_Lyapunov_and_S}), we calculate Lyapunov exponents for
variables $X_{i}$ such that $X_{i}^{\ast }X_{i}$ are distributed as Free
Poisson variables with parameter $\lambda .$ The case $\lambda =1$
corresponds to the random matrix case considered by C. M. Newman in %
\citep{newman86a}, and the results of this paper are in agreement with
Newman's \textquotedblleft triangle\textquotedblright\ law. In addition, our
results regarding the largest Lyapunov exponent agree with the results
regarding the norm of products of large random matrices in %
\citep{cohen_newman84}. Finally, our formula for computation of Lyapunov
exponents seems to be easier to apply than the non-linear integral
transformation developed in \citep{newman86a}.

An interesting by-product of our results is a relation between the extended
Fuglede-Kadison determinant and the Voiculescu $S$-transform, which allows
expressing each of them in terms of the other. In particular, if $Y$ is a
positive operator and if $\left\{ P_{t}\right\} $ is a family of projections
which are free of $Y$ and such that $E\left( P_{t}\right) =t$, then 
\begin{equation}
\log S_{Y}\left( -t\right) =-2\frac{d}{dt}\left[ \log \det \left( \sqrt{Y}%
P_{t}\right) \right] .  \label{formula_S_and_det}
\end{equation}%
and if $X$ is bounded and invertible, then 
\begin{equation}
\log \det \left( X\right) =-\frac{1}{2}\int_{0}^{1}\log S_{X^{\ast }X}\left(
-t\right) dt.  \label{formula_det_and_S}
\end{equation}

Calculations related to (\ref{formula_S_and_det}) and (\ref%
{formula_det_and_S}) were performed by Haagerup and Larsen in %
\citep{haagerup_larsen00} in their investigation of the Brown measure of $R$%
-diagonal operators. The Brown measure of an operator $X$ is closely related
to the determinant of $X-zI,$ and Haagerup and Larsen computed the Brown
measure of an $R$-diagonal operator $X$ in terms of the $S$-transform of $%
X^{\ast }X.$ However, it appears that formulas (\ref{formula_S_and_det}) and
(\ref{formula_det_and_S}) have not been stated explicitly in %
\citep{haagerup_larsen00}.

In addition, Sniady and Speicher showed in \citep{sniady_speicher01} that an 
$R$-diagonal $X$ can be represented in the triangular form and that the
spectra of the diagonal elements in this representation satisfy certain
inequalities in terms of the $S$-transform of $X^{\ast }X.$ Sniady and
Speicher used their result to give a different proof of the results in %
\citep{haagerup_larsen00}. It is likely that Sniady and Speicher's method
can also be used for a different proof of formulas (\ref{formula_S_and_det})
and (\ref{formula_det_and_S}).

The rest of the paper is organized as follows:\ Section \ref%
{section_FK_determinant} describes the extension of the Fuglede-Kadison
determinant that we use in this paper. Section \ref%
{section_existence_of_Lyapunov_exponents} defines the Lyapunov exponents of
free operators, proves an existence theorem, and derives a formula for the
calculation of Lyapunov exponents. Section \ref{section_example} computes
the Lyapunov exponents for a particular example. Section \ref%
{section_connection_Lyapunov_Stransform} connects the marginal Lyapunov
exponents and the $S$-transform, proves additivity and monotonicity of the
marginal Lyapunov exponent, and derives a formula for the largest Lyapunov
exponent. In addition, it derives a relation between the determinant and the 
$S$-transform. And Section \ref{section_conclusion} concludes.

\section{A modification of the Fuglede-Kadison determinant\label%
{section_FK_determinant}}

Let $\mathcal{A}$ be a finite von Neumann algebra and $E$ be a trace in this
algebra. Recall that if $X$ is an element of $\mathcal{A}$ that has a
bounded inverse, then the Fuglede-Kadison determinant (\cite%
{fuglede_kadison52}) is defined by the following formula:%
\begin{equation}
\det \left( X\right) =\exp \frac{1}{2}E\log \left( X^{\ast }X\right) .
\label{formula_FK_determinant_def1}
\end{equation}%
The most important property of the Fuglede-Kadison determinant is its
multiplicativity:%
\begin{equation}
\det \left( XY\right) =\det \left( X\right) \det \left( Y\right) .
\label{formula_FK_determinant_multiplicativity}
\end{equation}%
This determinant cannot be extended (non-trivially) to operators with
non-zero kernel if we require that property (\ref%
{formula_FK_determinant_multiplicativity}) holds for all $X$ and $Y.$

However, if we do not insist on this property, then we can define an \emph{%
extended determinant} as follows: Let $\log ^{+\lambda }\left( t\right)
=:\log t$ if $t>\lambda $ and $=:0$ if $t\leq \lambda .$ Note that $E\log
^{+\lambda }\left( X^{\ast }X\right) $ is a (weakly) decreasing function of $%
\lambda $ on the interval $\left( 0,1\right) ,$ and therefore it converges
to a limit (possibly infinite) as $\lambda \rightarrow 0.$

\begin{definition}
\label{definition_determinant}%
\begin{equation*}
\det \left( X\right) =\exp \frac{1}{2}\lim_{\lambda \downarrow 0}E\log
^{+\lambda }\left( X^{\ast }X\right) .
\end{equation*}
\end{definition}

This extension of the Fuglede-Kadison determinant coincides with the
extension introduced in Section 3.2 of \citep{luck02} .

\begin{example}
Zero operator
\end{example}

From Definition \ref{definition_determinant}, if $X=0,$ then $\det X=1.$

\begin{example}
Finite dimensional operator
\end{example}

Consider the algebra of $n$-by-$n$ matrices $M_{n}\left( C\right) $ with the
trace given as the normalization of the usual matrix trace: $E\left(
X\right) =n^{-1}\limfunc{Tr}\left( X\right) .$ Then the original
Fuglede-Kadison determinant is defined for all full-rank matrices and equals
the product of the singular values of the operator in the power of $1/n$. It
is easy to see that this equals the absolute value of the usual matrix
determinant in the power of $1/n$. The extended Fuglede-Kadison determinant
is defined for all matrices, including the matrices of rank $k<n,$ and
equals the product of non-zero singular values in the power of $1/n.$

\bigskip

We can write the definition of the determinant in a slightly different form.
Recall that for a self-adjoint operator $X\in \mathcal{A}$ we can define its 
\emph{spectral probability measure }as follows: 
\begin{equation*}
\mu _{X}\left( S\right) =E\left( 1_{S}\left( X\right) \right) ,
\end{equation*}%
where $S$ is an arbitrary Borel-measurable set and $1_{S}$ is its indicator
function. Then, the determinant of operator $X$ can be written as 
\begin{equation*}
\det \left( X\right) =\exp \frac{1}{2}\lim_{\lambda \downarrow 0}\int_{%
\mathbb{R}^{+}}\log ^{+\lambda }\left( t\right) \mu _{X^{\ast }X}\left(
dt\right) .
\end{equation*}

For all invertible $X$ the extended determinant defines the same object as
the usual Fuglede-Kadison determinant. For non-invertible $X$, the
multiplicativity property sometimes fails. However, it holds if a certain
condition on images and domains of the multiplicands is fulfilled:

\begin{proposition}
\label{theorem_determinant_multiplicativity}Let $V$ be the closure of the
range of the operator $X.$ If $Y$ is an injective mapping on $V$ and is the
zero operator on $V^{\perp }$, then $\det \left( YX\right) =\det \left(
Y\right) \det \left( X\right) .$
\end{proposition}

The claim of this proposition is a direct consequence of Theorem 3.14 and
Lemma 3.15(7) in \citep{luck02}.

\bigskip

\section{Definition of Lyapunov exponents for free operators\label%
{section_existence_of_Lyapunov_exponents}}

A pair $\left( \mathcal{A},E\right) $ is a \emph{tracial} $W^{\ast }$-\emph{%
non-commutative probability space} if $\mathcal{A}$ is a finite von Neumann
algebra with a normal faithful tracial state $E$, and $E\left( I\right) =1.$
The trace $E$ will be called the \emph{expectation} by analogy with
classical probability theory.

Let $\mathcal{A}_{1,}...,\mathcal{A}_{n}$ be sub-algebras of algebra $%
\mathcal{A}$, and let $a_{i}$ be elements of these sub-algebras such that $%
a_{i}\in \mathcal{A}_{k\left( i\right) }.$

\begin{definition}
The sub-algebras $\mathcal{A}_{1,}...,\mathcal{A}_{n}$ (and their elements)
are called \emph{free} or \emph{freely independent} if $E\left(
a_{1}...a_{m}\right) =0$ whenever the following two conditions hold:\newline
(a) $E\left( a_{i}\right) =0$ for every $i$, and \newline
(b) $k(i)\neq k\left( i+1\right) $ for every $i<m.$
\end{definition}

The random variables are called free or freely independent if the algebras
that they generate are free. (See \citep{voiculescu_dykema_nica92}, %
\citep{hiai_petz00}, or \citep{nica_speicher06} for more details on
foundations of free probability theory.)

Let $\left\{ X_{i}\right\} _{i=1}^{\infty }$ be a sequence of free
identically-distributed operators. Let $\Pi _{n}=X_{n}\ldots X_{1},$ and let 
$P_{t}$ be a projection which is free of all $X_{i}$ and has the dimension $%
t,$ i.e., $E\left( P_{t}\right) =t.$

\begin{definition}
The \emph{integrated Lyapunov exponent} corresponding to the sequence $X_{i}$
is a real-valued function of $t\in \left[ 0,1\right] $ which is defined as
follows: 
\begin{equation*}
F\left( t\right) =\lim_{n\rightarrow \infty }\frac{1}{n}\log \det \left( \Pi
_{n}P_{t}\right) ,
\end{equation*}%
provided that the limit exists.
\end{definition}

\textbf{Remark:} In the case of random matrices, $\Pi _{n}$ is the product
of independent identically-distributed random matrices. In this case, it
turns out that the function defined analogously to $F\left( t\right) $
equals the sum of the $tN$ largest Lyapunov exponents divided by $N$, where $%
N$ is the dimension of the matrices and $t$ belongs to the set $\left\{
0/N,1/N,\ldots ,N/N\right\} .$

Our first task is to prove the existence of the limit in the previous
definition.

\begin{theorem}
\label{theorem_Lyapunov_function}Suppose that $X_{i}$ are free
identically-distributed operators in a tracial $W^{\ast }$-probability space 
$\mathcal{A}$ with trace $E$. Let $u=:\dim \ker \left( X_{i}\right) .$ Then 
\begin{equation*}
F\left( t\right) =\left\{ 
\begin{array}{c}
\frac{1}{2}\lim_{\lambda \downarrow 0}E\log ^{+\lambda }\left(
P_{t}X_{1}^{\ast }X_{1}P_{t}\right) ,\text{ if }t\leq 1-u, \\ 
\frac{1}{2}\lim_{\lambda \downarrow 0}E\log ^{+\lambda }\left(
P_{u}X_{1}^{\ast }X_{1}P_{u}\right) ,\text{ if }t>1-u.%
\end{array}%
\right.
\end{equation*}
\end{theorem}

Before proving this theorem, let us make some remarks. First, this theorem
shows that the integrated Lyapunov exponent of the sequence $\left\{
X_{i}\right\} $ exists and depends only on the spectral distribution of $%
X_{i}^{\ast }X_{i}.$

Next, suppose that we know that $F\left( t\right) $ is differentiable almost
everywhere. Then we can define the \emph{marginal Lyapunov exponent} as $%
f\left( t\right) =F^{\prime }\left( t\right) .$ We can also define the \emph{%
distribution function of Lyapunov exponents }by the formula: $\mathcal{F}%
(x)=\mu \left\{ t\in \left[ 0,1\right] :f\left( t\right) \leq x\right\} ,$
where $\mu $ is the usual Borel-Lebesgue measure. Intuitively, this function
gives a measure of the set of the Lyapunov exponents which are less than a
given threshold, $x.$ In the finite-dimensional case it is simply the
empirical distribution function of the Lyapunov exponents, i.e., the
fraction of Lyapunov exponents that fall below the threshold $x.$

\textbf{Proof of Theorem \ref{theorem_Lyapunov_function}:} The proof is
through a sequence of lemmas. We will consider first the case of injective
operators $X_{i}$ and then will show how to generalize the argument to the
case of arbitrary $X_{i}.$

Let $P_{A}$ denote the projection on the closure of the range of operator $%
A. $ In the following lemmas we always assume that operators belong to a
tracial $W^{\ast }$-probability space $\mathcal{A}$ with trace $E.$

\begin{lemma}
\label{Lemma_PAPt_and_Pt}Suppose that operator $A$ is injective, and that $%
P_{t}$ is a projection of dimension $t.$ Then projection $P_{AP_{t}}$ is
equivalent to $P_{t}$. In particular, $E\left( P_{AP_{t}}\right) =t.$
\end{lemma}

\textbf{Proof: }Recall that polar decomposition is possible in $\mathcal{A}$%
. (See Proposition II.3.14 on p. 77 in \citep{takesaki79} for details.)
Therefore, we can write $AP_{t}=WB,$ where $W$ is a partial isometry and $B$
is positive, and where both $W$ and $B$ belong to $\mathcal{A}$. By
definition, the range of $W$ is $\left[ \mathrm{Range}\left( AP_{t}\right) %
\right] ,$ and the domain of $W$ is $\left[ x:Bx=0\right] ^{\bot }=\left[
x:AP_{t}x=0\right] ^{\bot }=$ $\left[ x:P_{t}x=0\right] ^{\bot }=\left[ 
\mathrm{Range}\left( P_{t}\right) \right] .$ Therefore, $P_{AP_{t}}$ is
equivalent to $P_{t}$, with the equivalence given by the partial isometry $%
W. $ In particular, $\dim \left( P_{AP_{t}}\right) =\dim \left( P_{t}\right)
,$ i.e., $E\left( P_{AP_{t}}\right) =t.$ QED.

\begin{lemma}
\label{Lemma_PAPt_and_B_are_free}If $A,$ $A^{\ast },$ and $P_{t}$ are free
from an operator subalgebra $\mathcal{B},$ then $P_{AP_{t}}$ is free from $%
\mathcal{B}.$
\end{lemma}

\textbf{Proof:} $P_{AP_{t}}$ belongs to the $W^{\ast }$-algebra generated by 
$I,$ $A,$ $A^{\ast },$ and \thinspace $P_{t}.$ By assumption, this algebra
is free from $\mathcal{B}.$ Hence, $P_{AP_{t}}$ is also free from $\mathcal{B%
}.$ QED.

Let us use the notation $Q_{k}=P_{X_{k}...X_{1}P_{t}}$ for $k\geq 1$ and $%
Q_{0}=P_{t}$. Then by Lemma \ref{Lemma_PAPt_and_B_are_free}, $Q_{k}$ is free
from $X_{k+1}.$ Besides, if all $X_{i}$ are injective, then their product is
injective and, therefore, by Lemma \ref{Lemma_PAPt_and_Pt}, $Q_{k}$ is
equivalent to $P_{t}.$

\begin{lemma}
\label{lemma_determinant_of_product}If all $X_{i}$ are injective, then%
\begin{equation*}
\det \left( \Pi _{n}P_{t}\right) =\prod\nolimits_{i=1}^{n}\det \left(
X_{i}Q_{i-1}\right) .
\end{equation*}
\end{lemma}

\textbf{Proof}: Note that $\Pi _{n}P_{t}=X_{n}Q_{n-1}X_{n-1}\ldots
Q_{1}X_{1}Q_{0}.$ We will proceed by induction. We need only to prove that 
\begin{equation}
\det \left( X_{k+1}Q_{k}X_{k}\ldots Q_{1}X_{1}Q_{0}\right) =\det \left(
X_{k+1}Q_{k}\right) \det \left( X_{k}\ldots Q_{1}X_{1}Q_{0}\right) .
\label{formula_product_reduction}
\end{equation}%
Let $V_{k}$ be the closure of the range of $X_{k}\ldots Q_{1}X_{1}Q_{0}.$
Since $X_{k+1}$ is injective and $Q_{k}$ is the projector on $V_{k}$,
therefore $X_{k+1}Q_{k}$ is injective on $V_{k}$ and equal to zero on $%
V_{k}^{\perp }.$ Consequently, we can apply Proposition \ref%
{theorem_determinant_multiplicativity} and obtain (\ref%
{formula_product_reduction}). QED.

Now we are ready to prove Theorem \ref{theorem_Lyapunov_function} for the
case of injective $X_{i}$. Using Lemma \ref{lemma_determinant_of_product},
we write 
\begin{equation*}
n^{-1}\log \det \left( \Pi _{n}P_{t}\right) =\frac{1}{n}\sum%
\nolimits_{i=1}^{n}\log \det \left( X_{i}Q_{i-1}\right) .
\end{equation*}%
Note that $X_{i}$ are identically distributed by assumption, $Q_{i}$ have
the same dimension by Lemma \ref{Lemma_PAPt_and_Pt}, and $X_{i}$ and $%
Q_{i-1} $ are free by Lemma \ref{Lemma_PAPt_and_B_are_free}. This implies
that $\lim_{\lambda \downarrow 0}E\log ^{+\lambda }\left( Q_{i-1}X_{i}^{\ast
}X_{i}Q_{i-1}\right) $ does not depend on $i,$ and hence, $\det \left(
X_{i}Q_{i-1}\right) $ does not depend on $i.$ Hence, using $i=1$ we can
write:%
\begin{equation*}
n^{-1}\log \det \left( \Pi _{n}P_{t}\right) =\log \det \left(
X_{1}P_{t}\right) .
\end{equation*}

This finishes the proof for the case of injective $X_{i}.$ For the case of
non-injective $X_{i}$, i.e., for the case when $\dim \ker \left(
X_{i}\right) >0,$ we need the following lemma.

\begin{lemma}
\label{Lemma_kernel_product}Suppose that $P_{t}$ is a projection operator
free of $A$ and such that $E\left( P_{t}\right) =t.$ Then $\dim \ker \left(
AP_{t}\right) =\max \left\{ 1-t,\dim \ker \left( A\right) \right\} .$
\end{lemma}

\textbf{Proof}: Let $V=\left( \func{Ker}A\right) ^{\perp }$ and let $P_{V}$
be the projection on $V.$ Then $E\left( P_{V}\right) =1-\dim \func{Ker}A.$
Note that $Ax=0\Longleftrightarrow P_{V}x=0.$ Consequently, $%
AP_{t}x=0\Longleftrightarrow P_{V}P_{t}x=0.$ Therefore, we have:%
\begin{eqnarray*}
\dim \left\{ x:AP_{t}x=0\right\} &=&\dim \left\{ x:P_{V}P_{t}x=0\right\} \\
&=&\dim \left\{ x:P_{t}P_{V}P_{t}x=0\right\} .
\end{eqnarray*}

Since $P_{t}$ and $P_{V}$ are free, an explicit calculation of the
distribution of $P_{t}P_{V}P_{t}$ shows that 
\begin{equation*}
\dim \left\{ x:P_{t}P_{V}P_{t}x=0\right\} =\max \left\{ 1-t,1-\dim V\right\}
.
\end{equation*}%
QED.

Consider first the case when $0<\dim \ker X_{i}\leq 1-t.$ This case is very
similar to the case of injective $X_{i}.$ Using Lemma \ref%
{Lemma_kernel_product} we conclude that $\dim \func{Ker}\left(
X_{1}P_{t}\right) =1-t,$ and therefore that $E\left( P_{X_{1}P_{t}}\right)
=t.$ If, as before, we denote $P_{X_{1}P_{t}}$ as $Q_{1},$ then the
projection $Q_{1}$ is free from $X_{2},$ and $E\left( Q_{1}\right) =t.$

Similarly, we obtain that $E\left( P_{X_{2}Q_{1}}\right) =t.$ Proceeding
inductively, we define $Q_{k}=P_{X_{k}Q_{k-1}}$ and conclude that $Q_{k}$ is
free from $X_{k+1}$ and that $E\left( Q_{k}\right) =t.$

Next, we write $X_{k}...X_{1}P_{t}=X_{k}Q_{k-1}X_{k-1}Q_{k-2}...X_{1}Q_{0}$,
where $Q_{0}$ denotes $P_{t},$ and note that $X_{k}Q_{k-1}$ is injective on
the range of $Q_{k-1}.$ Indeed, if it were not injective, then we would have 
$\dim \left( \func{Ker}\left( X_{k}\right) \cap \mathrm{Range}\left(
Q_{k-1}\right) \right) >0.$ But this would imply that $\dim \left( \func{Ker}%
X_{k}Q_{k-1}\right) >\dim \left( \func{Ker}Q_{k-1}\right) =1-t,$ which
contradicts the fact that $\dim \left( \func{Ker}X_{k}Q_{k-1}\right) =1-t.$
Therefore, Proposition \ref{theorem_determinant_multiplicativity} is
applicable and 
\begin{eqnarray*}
\det \left( X_{k}...X_{1}P_{t}\right) &=&\det \left( X_{k}Q_{k-1}\right)
...\det \left( X_{1}Q_{0}\right) \\
&=&\left[ \det \left( X_{1}P_{t}\right) \right] ^{k}.
\end{eqnarray*}

Now let us turn to the case when $\dim \ker X_{i}=u>1-t.$ Then $\dim \func{%
Ker}\left( X_{1}P_{t}\right) =1-u$ and therefore $E\left(
P_{X_{1}P_{t}}\right) =u.$ Proceeding as before, we conclude that $E\left(
Q_{k}\right) =u$ for all $k\geq 1,$ and we can write $%
X_{k}...X_{1}P_{t}=X_{k}Q_{k-1}X_{k-1}Q_{k-2}...X_{1}Q_{0}$, where we have
denoted $P_{t}$ as $Q_{0}.$ Then we get the following formula:%
\begin{eqnarray*}
\det \left( X_{k}...X_{1}P_{t}\right) &=&\det \left( X_{k}Q_{k-1}\right)
...\det \left( X_{2}Q_{1}\right) \det \left( X_{1}Q_{0}\right) \\
&=&\left[ \det \left( X_{1}P_{u}\right) \right] ^{k-1}\det \left(
X_{1}P_{t}\right) .
\end{eqnarray*}

Therefore, 
\begin{equation*}
\lim_{n\rightarrow \infty }n^{-1}\log \det \left( \Pi _{n}P_{t}\right) =\log
\det \left( X_{1}P_{u}\right) .
\end{equation*}%
QED.

\section{Example \label{section_example}}

Let us compute the Lyapunov exponents for a random variable $X$ that has the
product $X^{\ast }X$ distributed according to the Marchenko-Pastur
distribution. Recall that the continuous part of the Marchenko-Pastur
probability distribution with parameter $\lambda >0$ is supported on the
interval $\left[ \left( 1-\sqrt{\lambda }\right) ^{2},\left( 1+\sqrt{\lambda 
}\right) ^{2}\right] ,$ and has the following density there:%
\begin{equation*}
p_{\lambda }\left( x\right) =\frac{\sqrt{4\lambda -\left( x-1-\lambda
\right) ^{2}}}{2\pi x}.
\end{equation*}%
For $\lambda \in \left( 0,1\right) ,$ this distribution also has an atom at $%
0$ with the probability mass $\left( 1-\lambda \right) $ assigned to it. The
Marchenko-Pastur distribution is sometimes called the free Poisson
distribution since it arises as a limit of free additive convolutions of the
Bernoulli distribution, \ and a similar limit in the classical case equals
the Poisson distribution. It can also be thought of as a scaled limit of the
eigenvalue distribution of Wishart-distributed random matrices (see \cite%
{hiai_petz00} for a discussion).

\begin{proposition}
\label{theorem_Lyapunov_for_MP}Suppose that $X$ is a non-commutative random
variable in a tracial $W^{\ast }$-probability space $\left( \mathcal{A}%
,E\right) ,$ such that $X^{\ast }X$ is distributed according to the
Marchenko-Pastur distribution with parameter $\lambda .$ If $\lambda \geq 1,$
then the distribution of Lyapunov exponents of $X$ is 
\begin{equation*}
\mathcal{F}\left( x\right) =\left\{ 
\begin{array}{cc}
0, & \text{ if }x<\left( 1/2\right) \log \left( \lambda -1\right) \\ 
e^{2x}+1-\lambda , & \text{if }x\in \left[ \frac{1}{2}\log \left( \lambda
-1\right) ,\frac{1}{2}\log \lambda \right) \\ 
1 & \text{if }x\geq \frac{1}{2}\log \lambda .%
\end{array}%
\right. .
\end{equation*}%
If $\lambda <1,$ then the distribution of Lyapunov exponents of $X$ is 
\begin{equation*}
\mathcal{F}\left( x\right) =\left\{ 
\begin{array}{cc}
e^{2x}, & \text{ if }x<\left( 1/2\right) \log \left( \lambda \right) \\ 
\lambda , & \text{if }x\in \left[ \frac{1}{2}\log \left( \lambda \right)
,0\right) \\ 
1 & \text{if }x\geq 0.%
\end{array}%
\right.
\end{equation*}
\end{proposition}

\textbf{Remark:} If $\lambda =1,$ then the distribution is the exponential
law discovered by C. M. Newman as a scaling limit of Lyapunov exponents of
large random matrices. (See \citep{newman86a}, \citep{newman86b}, and \cite%
{isopi_newman92}. This law is often called the \textquotedblleft
triangle\textquotedblright\ law since it implies that the exponentials of
Lyapunov exponents converge to the law whose density function is in the form
of a triangle.)

\textbf{Proof of Proposition \ref{theorem_Lyapunov_for_MP}: }It is easy to
calculate that the continuous part of the distribution of $P_{t}XP_{t}$ is
supported on the interval $\left[ \left( \sqrt{t}-\sqrt{\lambda }\right)
^{2},\left( \sqrt{t}+\sqrt{\lambda }\right) ^{2}\right] ,$ and has the
density function%
\begin{equation*}
p_{t,\lambda }\left( x\right) =\frac{\sqrt{4\lambda t-\left[ x-\left(
t+\lambda \right) \right] ^{2}}}{2\pi x}.
\end{equation*}%
This distribution also has an atom at $x=0$ with the probability mass $\max
\left\{ 1-\lambda ,1-t\right\} $. See for example, results in \cite%
{nica_speicher96}.

Next, we write the expression for the integrated Lyapunov exponent. If $%
\lambda \geq 1,$ or $\lambda <1$ but $\lambda \geq t,$ then 
\begin{eqnarray}
F_{\lambda }\left( t\right) &=&\frac{1}{2}\lim_{\varepsilon \downarrow
0}E\log ^{+\varepsilon }\left( P_{t}X^{\ast }XP_{t}\right)  \notag \\
&=&\frac{1}{2}\int\limits_{\left( \sqrt{t}-\sqrt{\lambda }\right)
^{2}}^{\left( \sqrt{t}+\sqrt{\lambda }\right) ^{2}}\log x\frac{\sqrt{%
4\lambda t-\left[ x-\left( t+\lambda \right) \right] ^{2}}}{2\pi x}dx.
\label{formula_MP_Lyapunov_distribution}
\end{eqnarray}%
If $\lambda <1$ and $\lambda <t,$ then 
\begin{eqnarray}
F_{\lambda }\left( t\right) &=&\frac{1}{2}\lim_{\varepsilon \downarrow
0}E\log ^{+\varepsilon }\left( P_{1-\lambda }X^{\ast }XP_{1-\lambda }\right)
\notag \\
F_{\lambda }\left( t\right) &=&\frac{1}{2}\int\limits_{\left( \sqrt{%
1-\lambda }-\sqrt{\lambda }\right) ^{2}}^{\left( \sqrt{1-\lambda }+\sqrt{%
\lambda }\right) ^{2}}\log x\frac{\sqrt{4\lambda \left( 1-\lambda \right) -%
\left[ x-1\right] ^{2}}}{2\pi x}dx.
\end{eqnarray}%
Differentiating (\ref{formula_MP_Lyapunov_distribution}) with respect to $t,$
we obtain an expression for the marginal Lyapunov exponent:%
\begin{equation}
f_{\lambda }\left( t\right) =\frac{1}{4\pi }\int\limits_{\left( \sqrt{t}-%
\sqrt{\lambda }\right) ^{2}}^{\left( \sqrt{t}+\sqrt{\lambda }\right) ^{2}}%
\frac{\log x}{x}\frac{x-t+\lambda }{\sqrt{4\lambda x-\left[ x-t+\lambda %
\right] ^{2}}}dx.  \label{formula_MP_Lyapunov_marginal}
\end{equation}

Using substitutions $u=\left[ x-\left( \sqrt{t}-\sqrt{\lambda }\right) ^{2}%
\right] /\left( 2\sqrt{\lambda t}\right) -1$ and then $\theta =\arccos u,$
this integral can be computed as 
\begin{equation*}
f_{\lambda }\left( t\right) =\frac{1}{2}\log \left( \lambda -t\right) .
\end{equation*}

From this expression, we calculate the distribution of Lyapunov exponents
for the case when $\lambda \geq 1$: 
\begin{equation*}
\mathcal{F}\left( x\right) =\left\{ 
\begin{array}{cc}
0, & \text{ if }x<\left( 1/2\right) \log \left( \lambda -1\right) \\ 
e^{2x}+1-\lambda , & \text{if }x\in \left[ \frac{1}{2}\log \left( \lambda
-1\right) ,\frac{1}{2}\log \lambda \right) \\ 
1 & \text{if }x\geq \frac{1}{2}\log \lambda .%
\end{array}%
\right.
\end{equation*}

A similar analysis shows that for $\lambda <1,$ the distribution is as
follows:%
\begin{equation*}
\mathcal{F}\left( x\right) =\left\{ 
\begin{array}{cc}
e^{2x}, & \text{ if }x<\left( 1/2\right) \log \left( \lambda \right) \\ 
\lambda , & \text{if }x\in \left[ \frac{1}{2}\log \left( \lambda \right)
,0\right) \\ 
1 & \text{if }x\geq 0.%
\end{array}%
\right.
\end{equation*}%
QED.

\section{A relation with the $S$-transform \label%
{section_connection_Lyapunov_Stransform}}

In this section we derive a formula that makes the calculation of Lyapunov
exponents easier and relates them to the $S$-transform of the operator $%
X_{i} $. Recall that the $\psi $-function of a bounded non-negative operator 
$A$ is defined as $\psi _{A}\left( z\right) =\sum_{k=1}^{\infty }E\left(
A^{k}\right) z^{k}.$ Then the $S$-transform is $S_{A}\left( z\right) =\left(
1+z^{-1}\right) \psi _{A}^{\left( -1\right) }\left( z\right) ,$ where $\psi
_{A}^{\left( -1\right) }\left( z\right) $ is the functional inverse of $\psi
_{A}\left( z\right) $ in a neighborhood of $0.$

\begin{theorem}
\label{theorem_marginal_Lyapunov_function}Let $X_{i}$ be identically
distributed free bounded operators in a tracial $W^{\ast }$-probability
space $\mathcal{A}$ with trace $E$. Let $Y=X_{1}^{\ast }X_{1}$ and suppose
that $\mu _{Y}\left( \left\{ 0\right\} \right) =1-u\in \left[ 0,1\right) ,$
where $\mu _{Y}$ denotes the spectral probability measure of $Y.$ Then the
marginal Lyapunov exponent of the sequence $\left\{ X_{i}\right\} $ is given
by the following formula:%
\begin{equation*}
f_{X}\left( t\right) =\left\{ 
\begin{array}{cc}
-\frac{1}{2}\log \left[ S_{Y}\left( -t\right) \right] & \text{if }t<u, \\ 
0 & \text{if }t>u,%
\end{array}%
\right.
\end{equation*}%
where $S_{Y}$ is the $S$-transform of the variable $Y.$
\end{theorem}

\textbf{Remark:} Note that if $X_{1}^{\ast }X_{1}$ has no atom at zero then
the formula is simply $f_{X}\left( t\right) =-\frac{1}{2}\log \left[
S_{Y}\left( -t\right) \right] .$

\textbf{Proof:} If $t>u,$ then $f_{X}\left( t\right) =0$ by Theorem \ref%
{theorem_Lyapunov_function}. Assume in the following that $t<u.$ Then $%
P_{t}X^{\ast }XP_{t}$ has an atom of mass $1-t$ at $0.$ Let $\mu _{t}$
denote the spectral probability measure of $P_{t}X^{\ast }XP_{t},$ with the
atom at $0$ removed. (So the total mass of $\mu _{t}$ is $t$.) We start with
the formula:%
\begin{equation*}
\log x=\log \left( c+x\right) -\int_{0}^{c}\frac{ds}{x+s},
\end{equation*}%
and write:%
\begin{equation*}
\int_{0}^{\infty }\log x\text{ }\mu _{t}\left( dx\right) =\lim_{c\rightarrow
\infty }\left\{ t\log \left( c\right) +\int_{0}^{c}G_{t}\left( -s\right)
ds\right\} ,
\end{equation*}%
where $G_{t}$ is the Cauchy transform of the measure $\mu _{t}.$

Next, note that $G_{t}\left( -s\right) =-s^{-1}\psi _{t}\left(
-s^{-1}\right) -ts^{-1}$ and substitute this into the previous equation: 
\begin{eqnarray*}
\int_{0}^{\infty }\log x\text{ }\mu _{t}\left( dx\right)
&=&\lim_{c\rightarrow \infty ,\varepsilon \rightarrow 0}\left\{ t\log
c-t\log c+t\log \left( \varepsilon \right) +\int_{\varepsilon }^{c}\frac{%
\psi _{t}\left( s^{-1}\right) }{s}ds\right\} \\
&=&\lim_{\varepsilon \rightarrow 0}\left\{ t\log \left( \varepsilon \right)
+\int_{\varepsilon }^{\infty }\frac{\psi _{t}\left( s^{-1}\right) }{s}%
ds\right\} .
\end{eqnarray*}%
Using substitutions $v=-\log s,$ and $A=-\log \varepsilon ,$ we can re-write
this equation as follows:%
\begin{equation*}
\int_{0}^{\infty }\log x\text{ }\mu _{t}\left( dx\right) =\lim_{A\rightarrow
\infty }\left\{ -tA-\int_{-\infty }^{A}\psi _{t}\left( -e^{v}\right)
dv\right\} .
\end{equation*}

The function $\psi _{t}\left( -e^{v}\right) $ monotonically decreases when $%
v $ changes from $-\infty $ to $\infty ,$ and its value changes from $0$ to $%
-t.$ Let $s^{\ast }=:\psi _{t}\left( -e^{0}\right) =\psi _{t}\left(
-1\right) $ and let $\xi _{t}\left( x\right) $ denote the functional inverse
of $\psi _{t}\left( -e^{v}\right) .$ The function $\xi _{t}\left( x\right) $
is defined on the interval $\left( -t,0\right) .$ In this interval it is
monotonically decreasing from $\infty $ to $-\infty .$ The only zero of $\xi
_{t}\left( x\right) $ is at $x=s^{\ast }.$

It is easy to see that 
\begin{equation*}
\lim_{A\rightarrow \infty }\left\{ -tA-\int_{0}^{A}\psi _{t}\left(
-e^{v}\right) dv\right\} =-\int_{-t}^{s^{\ast }}\xi _{t}\left( x\right) dx,
\end{equation*}%
and that 
\begin{equation*}
-\int_{-\infty }^{0}\psi _{t}\left( -e^{v}\right) dv=-\int_{s^{\ast
}}^{0}\xi _{t}\left( x\right) dx.
\end{equation*}%
Therefore, 
\begin{equation*}
\int_{0}^{\infty }\log x\text{ }\mu _{t}\left( dx\right) =-\int_{-t}^{0}\xi
_{t}\left( x\right) dx.
\end{equation*}%
It remains to note that $\xi _{t}\left( x\right) =\log \left[ -\psi
_{t}^{\left( -1\right) }\left( x\right) \right] ,$ in order to conclude that 
\begin{equation*}
\int_{0}^{\infty }\log x\text{ }\mu _{t}\left( dx\right) =-\int_{-t}^{0}\log %
\left[ -\psi _{t}^{\left( -1\right) }\left( x\right) \right] dx.
\end{equation*}%
The next step is to use Voiculescu's multiplication theorem and write: $\psi
_{t}^{\left( -1\right) }\left( x\right) =\psi _{Y}^{\left( -1\right) }\left(
x\right) \left( 1+x\right) /\left( t+x\right) .$ Then we have the formula:%
\begin{eqnarray*}
\int_{0}^{\infty }\log x\text{ }\mu _{t}\left( dx\right)
&=&-\int_{-t}^{0}\log \left[ -\psi _{Y}^{\left( -1\right) }\left( x\right) %
\right] dx-\int_{-t}^{0}\log \left[ \frac{1+x}{t+x}\right] dx \\
&=&-\int_{-t}^{0}\log \left[ -\psi _{Y}^{\left( -1\right) }\left( x\right) %
\right] dx+\left( 1-t\right) \log \left( 1-t\right) +t\log t.
\end{eqnarray*}%
The integrated Lyapunov exponent is one half of this expression, and we can
obtain the marginal Lyapunov exponent by differentiating over $t$:%
\begin{eqnarray*}
f\left( t\right) &=&\frac{1}{2}\left( -\log \left[ -\psi _{Y}^{\left(
-1\right) }\left( -t\right) \right] +\log t-\log \left( 1-t\right) \right) \\
&=&-\frac{1}{2}\log \left[ \left( 1-\frac{1}{t}\right) \psi _{Y}^{\left(
-1\right) }\left( -t\right) \right] \\
&=&-\frac{1}{2}\log \left[ S_{Y}\left( -t\right) \right] .
\end{eqnarray*}

QED.

\textbf{Example}

Let us consider again the case of identically distributed free $X_{i}$ such
that $X_{i}^{\ast }X_{i}$ has the Marchenko-Pastur distribution with the
parameter $\lambda \geq 1.$ In this case $S_{Y}\left( z\right) =\left(
\lambda +z\right) ^{-1}.$ Hence, applying Theorem \ref%
{theorem_marginal_Lyapunov_function}, we immediately obtain a formula for
the marginal Lyapunov exponent:%
\begin{equation*}
f\left( t\right) =\frac{1}{2}\log \left( \lambda -t\right) .
\end{equation*}%
Inverting this formula, we obtain the formula for the distribution of
Lyapunov exponents:%
\begin{equation*}
\mathcal{F}\left( x\right) =\left\{ 
\begin{array}{cc}
0, & \text{ if }x<\left( 1/2\right) \log \left( \lambda -1\right) , \\ 
e^{2x}+1-\lambda , & \text{if }x\in \left[ \frac{1}{2}\log \left( \lambda
-1\right) ,\frac{1}{2}\log \lambda \right) , \\ 
1 & \text{if }x\geq \frac{1}{2}\log \lambda ,%
\end{array}%
\right.
\end{equation*}%
which is exactly the formula that we obtained earlier by a direct
calculation from definitions. It is easy to check that a similar agreement
holds also for $\lambda <1.$

In the following corollaries we always assume that operators belong to a
tracial $W^{\ast }$-probability space $\mathcal{A}$ with trace $E.$

\begin{corollary}
Let $X$ and $Y$ be free and such that $X^{\ast }X$ and $Y^{\ast }Y$ are
bounded and have no atom at zero. Let $f_{X},$ $f_{Y},$ and $f_{XY}$ denote
the marginal Lyapunov exponents corresponding to variables $X,$ $Y$ and $XY,$
respectively. Then 
\begin{equation*}
f_{XY}\left( t\right) =f_{X}\left( t\right) +f_{Y}\left( t\right) .
\end{equation*}
\end{corollary}

\textbf{Proof:} By Theorem \ref{theorem_marginal_Lyapunov_function},%
\begin{eqnarray*}
f_{XY}\left( t\right) &=&-\frac{1}{2}\log \left[ S_{Y^{\ast }X^{\ast
}XY}\left( -t\right) \right] \\
&=&-\frac{1}{2}\log \left[ S_{Y^{\ast }Y}\left( -t\right) S_{X^{\ast
}X}\left( -t\right) \right] \\
&=&f_{X}\left( t\right) +f_{Y}\left( t\right) .
\end{eqnarray*}

QED.

\begin{corollary}
If $X$ is bounded and $X^{\ast }X$ has no atom at zero, then the marginal
Lyapunov exponent is (weakly) decreasing in $t$, i.e. $f_{X}^{\prime }\left(
t\right) \leq 0.$
\end{corollary}

\textbf{Proof:} Because of Theorem \ref{theorem_marginal_Lyapunov_function},
we need only to check that $S\left( t\right) $ is (weakly) decreasing on the
interval $\left[ -1,0\right] ,$ and this was proved by Bercovici and
Voiculescu in Proposition 3.1 on page 225 of \citep{bercovici_voiculescu92}.
QED.

\begin{corollary}
\label{theorem_leading_Lyapunov_exponent}If $X$ is bounded and $X^{\ast }X$
has no atom at zero, then the largest Lyapunov exponent equals $\left(
1/2\right) \log E(X^{\ast }X).$
\end{corollary}

\textbf{Proof:} This follows from the previous Corollary and the fact that $%
S_{Y}\left( 0\right) =1/E\left( Y\right) .$ QED.

\textbf{Remark: }It is interesting to compare this result with the result in %
\citep{cohen_newman84}, which shows that the norm of the product of $N\times
N$ i.i.d. random matrices $X_{1},\ldots ,X_{n}$ grows exponentially when $n$
increases, and that the asymptotic growth rate approaches $\frac{1}{2}\log E(%
\mathrm{tr}\left( X_{1}^{\ast }X_{1}\right) )$ if $N\rightarrow \infty $ and
matrices are scaled appropriately. The assumption in \citep{cohen_newman84}
about the distribution of matrix entries is that the distribution of $%
X_{1}^{\ast }X_{1}$ is invariant relative to orthogonal rotations of the
ambient space. Since the growth rate of the norm of the product $X_{1}\ldots
X_{n}$ is another way to define the largest Lyapunov exponent of the
sequence $X_{i},$ therefore the result in \citep{cohen_newman84} is in
agreement with Corollary \ref{theorem_leading_Lyapunov_exponent}.

The main result of Theorem \ref{theorem_marginal_Lyapunov_function} can also
be reformulated as the following interesting identity:

\begin{corollary}
\label{theorem_leading_Lyapunov_exponent copy(1)}If $Y$ is bounded,
self-adjoint, and positive, and if $\left\{ P_{t}\right\} $ is a family of
projections which are free of $Y$ and such that $E\left( P_{t}\right) =t$,
then 
\begin{eqnarray*}
\log S_{Y}\left( -t\right) &=&-\frac{d}{dt}\left[ \lim_{\lambda \rightarrow
0}E\log ^{+\lambda }\left( P_{t}YP_{t}\right) \right] \\
&=&-2\frac{d}{dt}\left[ \log \det \left( \sqrt{Y}P_{t}\right) \right] .
\end{eqnarray*}
\end{corollary}

Conversely, we can express the determinant in terms of the $S$-transform:

\begin{corollary}
If $X$ is bounded and invertible, then 
\begin{equation*}
\det \left( X\right) =\exp \left\{ -\frac{1}{2}\int_{0}^{1}\log S_{X^{\ast
}X}\left( -t\right) dt\right\} .
\end{equation*}
\end{corollary}

\section{Conclusion \label{section_conclusion}}

In conclusion we want to indicate how our results are related to results by
Newman in \citep{newman86a}. Newman considers $N$-by-$N$ random matrices $%
X_{i}$ with rotationally-invariant distribution.of entries. Suppose that the
empirical distribution of eigenvalues of $\sqrt{X_{i}^{\ast }X_{i}}$
converges to a probability measure $K\left( dx\right) $ as $N\rightarrow
\infty .$ Let $H\left( x\right) $ be the limit cumulative distribution
function for $e^{f_{k}},$ where $f_{k}$ are Lyapunov exponents of $X_{i}.$
Newman shows that this limit exists and satisfies the following integral
equation: 
\begin{equation*}
\int\limits_{0}^{\infty }\frac{t^{2}}{H\left( x\right) x^{2}+\left(
1-H\left( x\right) \right) t^{2}}K\left( dt\right) =1.
\end{equation*}

We are going to transform this expression in a form which is more comparable
with our results. Let $\mu \left( dt\right) $ be the limit of the empirical
distribution of eigenvalues of $X_{i}^{\ast }X_{i}$ as $N\rightarrow \infty
, $ and let $s=t^{2}.$ Then, we can re-write Newman's formula as follows:%
\begin{equation*}
\int\limits_{0}^{\infty }\frac{s}{H\left( x\right) x^{2}+\left( 1-H\left(
x\right) \right) s}\mu \left( ds\right) =1.
\end{equation*}%
This expression can be transformed to the following equality: 
\begin{equation*}
\psi _{Y}\left( \frac{H\left( x\right) -1}{H\left( x\right) x^{2}}\right)
=H\left( x\right) -1,
\end{equation*}%
where $\psi _{Y}\left( z\right) $ is the $\psi $-function of a
non-commutative random variable $Y$ with the spectral probability
distribution $\mu $. Hence 
\begin{equation*}
\frac{1}{x^{2}}=S_{Y}\left[ H\left( x\right) -1\right] .
\end{equation*}

If we denote $1-H\left( x\right) $ as $u,$ then 
\begin{equation*}
x=e^{-\frac{1}{2}\log \left[ S_{Y}\left( -u\right) \right] }.
\end{equation*}%
This implies that as $N$ grows, the ordered Lyapunov exponents $f_{k}$
converge to $-\frac{1}{2}\log \left[ S_{Y}\left( -u\right) \right] ,$
provided that $k/N\rightarrow u.$ This is in agreement with our formula in
Theorem \ref{theorem_marginal_Lyapunov_function}.

\bibliographystyle{abbrvnat}
\bibliography{comtest}

\end{document}